\documentclass[12pt]{article}

\usepackage{tikz,color,url}
\usepackage{amsmath,amssymb,latexsym}

\usepackage{algorithm}
\usepackage{algpseudocode}
\algnotext{EndFor}
\algnotext{EndIf}

\newtheorem{theorem}{Theorem}[section]

\newtheorem{lemma}[theorem]{Lemma}

\newcommand{\proof}{\noindent{\bf Proof.\ }}
\newcommand{\qed}{\hfill $\square$ \medskip}
\newcommand{\cp}{\,\square\,}

\newcommand{\dimm}{{\rm dim}_{\rm ms}}
\newcommand{\diam}{{\rm diam}}
\newcommand{\ecc}{{\rm ecc}}
\newcommand{\Tr}{{\rm Tr}}

\newcommand{\msl}{\{\hspace*{-0.1cm}\{}
\newcommand{\msr}{\}\hspace*{-0.1cm}\}}

\begin{document}

\title{Further contributions on the outer multiset dimension of graphs}

\author{
Sandi Klav\v zar$^{,a,b,c}$
\and
Dorota Kuziak $^d$
\and
Ismael G. Yero $^e$
}

\date{}

\maketitle

\begin{center}
$^a$ Faculty of Mathematics and Physics, University of Ljubljana, Slovenia \\
\medskip

$^b$ Faculty of Natural Sciences and Mathematics, University of Maribor, Slovenia \\
\medskip

$^{c}$ Institute of Mathematics, Physics and Mechanics, Ljubljana, Slovenia \\
{\tt sandi.klavzar@fmf.uni-lj.si}
\medskip

$^{d}$ Departamento de Estad\'istica e Investigaci\'on Operativa, Universidad de C\'adiz, Algeciras, Spain \\
{\tt dorota.kuziak@uca.es}
\medskip

$^{e}$ Departamento de Matem\'aticas, Universidad de C\'adiz, Algeciras, Spain \\
{\tt ismael.gonzalez@uca.es}
\end{center}

\begin{abstract}
The outer multiset dimension ${\rm dim}_{\rm ms}(G)$ of a graph $G$ is the cardinality of a smallest set of vertices that uniquely recognize all the vertices outside this set by using multisets of distances to the set. It is proved that ${\rm dim}_{\rm ms}(G) = n(G) - 1$ if and only if $G$ is a regular graph with diameter at most $2$. Graphs $G$ with ${\rm dim}_{\rm ms}(G)=2$ are described and recognized in polynomial time. A lower bound on the lexicographic product of $G$ and $H$ is proved when $H$ is complete or edgeless, and the extremal graphs are determined. It is proved that ${\rm dim}_{\rm ms}(P_s\,\square\, P_t) = 3$ for $s\ge t\ge 2$.
\end{abstract}


\noindent
{\bf AMS Subj.\ Class.\ (2020)}: 05C12; 05C76

\section{Introduction}
\label{sec:intro}

Computing the metric dimension of graphs is one of the classical topics in metric graph theory, based on its applicability in several location related problems arising in different areas of investigation like for instance computer science, chemistry, biology, or social sciences. For a better understanding on definitions, terminology, contributions, and open questions on this issue we suggest the fairly complete and recently presented survey \cite{till-2022+}.

The theory of metric dimension in graphs has also been much developed by studying several variations of the classical concept, while attempting to consider some specifications or generalizations that are giving more insight into this classical version. The number of such variants has been significantly increased in the last recent years, and the reader can now find a very rich area of investigation concerning such variants. For a comprehensive background on a large number of such variants, it is suggested the other recent survey \cite{dorota-2022+}.

Given a connected graph $G$, it is said that a vertex $v\in V(G)$ \emph{resolves} (or \emph{identifies}, or \emph{determines}) two vertices $x,y\in V(G)$ if $d_G(x,v)\ne d_G(y,v)$; equivalently, $x,y$ are {\em resolved} by $v$. Here and later $d_G(u,w)$ stands for the distance between $u,w$ in $G$. It is also said that a set of vertices $S$ \emph{resolves} the set $V(G)$ if every two vertices of $G$ are resolved by a vertex of $S$, and such set is called a \emph{resolving set}. The \emph{metric dimension} of $G$ is defined as the cardinality of a smallest resolving set for $G$, and denoted by $\dim(G)$. A resolving set of cardinality $\dim(G)$ is called a \emph{metric basis}. The concepts above were first (and independently) presented in \cite{Harary,Slater}.

It can be easily noted that a resolving set $S$ of a graph $G$ has the property of uniquely identifying all the vertices of $G$ by means of distances to the vertices in $S$. That is, consider $S=\{v_1,\dots,v_k\}$ as an ordered set of vertices of a connected graph $G$, and for any vertex $v\in V(G)$, consider the vector
$$r(v|S)=(d_G(v,v_1),\dots,d_G(v,v_k)).$$
It is readily seen that $S$ is a resolving set for $G$ if and only if all the vectors $r(v|S)$ with $v\in V(G)$ are pairwise different. The vector $r(v|S)$ is called the \emph{metric representation} of $v$ with respect to $S$.

The core of the locating property of a resolving set $S$ of a graph is based on the uniqueness of the metric representations with respect to $S$ of the vertices of the graph, and these metric representations are given by vectors of distances. A modified version of this was presented in \cite{rino-2017}, where the authors suggested the use of ``multisets'' instead of vectors in the definition of metric representations of a vertex with respect to a given set. That is, for a given vertex $v\in V(G)$ and a set $S = \{w_1, \ldots, w_t\}$, the {\em multiset
representation of $u$ with respect to} $S$ is
$${\rm m}_G(u|S)=\msl d_G(u, w_1), \ldots, d_G(u, w_t) \msr,$$
where $\msl \cdot \msr$ limits a multiset. Hence, the set $S$ is a \emph{multiset resolving set} for $G$ if all the multisets ${\rm m}_G(u|S)$ with $u\in V(G)$ are pairwise different. The \emph{multiset dimension} of $G$ is then defined as the cardinality of a smallest multiset resolving set. A detail that one can immediately notice is that there could be vertices in a graph that have the same multiset representation with respect to every set of vertices of $G$ (for instance twin vertices), and so, the multiset dimension of such graphs is not properly defined. In such situations, the authors of \cite{rino-2017} adopted the agreement that such graphs has infinite multiset dimension. The problem from~\cite{rino-2017} to characterize the graphs with infinite multiset dimension remains open. Some partial contributions on this direction were already described in \cite{bong-2021}.

On the other hand, in order to avoid the problem of the possible infiniteness of the multiset version of the metric dimension, it was introduced an ``outer'' version of multiset resolving sets in \cite{gil-pons-2019} as follows. A set $S\subseteq V(G)$ is a \emph{outer multiset resolving set} for $G$, if the multiset representations of vertices $u\notin S$ with respect to $S$ are pairwise different. A multiset resolving set of the smallest possible cardinality is called an \emph{outer multiset basis}, and the cardinality of an outer multiset basis is the \emph{outer multiset dimension} of $G$, denoted by $\dimm(G)$. This structure clearly avoids the problem of infiniteness of multiset dimension, since vertices that must be distinguished by a set $S$ are only vertices outside $S$. In this work, we are aimed to continue developing this research line. To this end, we next give some basic definitions and terminologies that are necessary along our exposition.

For a positive integer $k$ we will use the notation $[k] = \{1,\ldots, k\}$. A vertex $u$ of a graph $G$ is {\em diametral} if there exists a vertex $v\in V(G)$ such that $d_G(u,v) = \diam(G)$, where $\diam(G)$ denotes the {\em diameter} of $G$, that is, the largest distance between vertices of $G$. We also say that $v$ is a vertex diametral to $u$. A subgraph $H$ of a graph $G$ is {\em isometric} if $d_H(u,v) = d_G(u,v)$ holds for all $u,v\in V(H)$. The open and the closed neighborhood of a vertex $u$ of $G$ will be denoted by $N_G(u)$ and by $N_G[u]$, respectively. The {\em degree} $\deg_G(u)$ of $u$ is $|N_G(u)|$. Vertices $x$ and $y$ of $G$ are {\em true twins} if $N_G[u] = N_G[v]$ and are {\em false twins} if $N_G(u) = N_G(v)$. Vertices $u$ and $v$ are just {\em twins}, if they are either true twins or false twins. The order of $G$ will be denoted by $n(G)$. If in a multiset (of distances) an element $d$ appears $k$ times, then we may abbreviate it to $d^k$. For instance, $\msl 1, 1, 1, 2, 5, 5 \msr = \msl 1^3, 2, 5^2 \msr$. Finally, all graphs considered in this paper are connected and of order at least $2$.

\section{Graphs with outer multiset dimension order minus one}
\label{sec:minus-one}

For the classical metric dimension parameter, it is well known that $\dim(G)=n(G)-1$ if and only if $G$ is a complete graphs. In \cite{gil-pons-2019}, some examples of non complete graphs $G$ of order $n$ such that $\dimm(G)=n-1$ were given. However, a complete characterization of the class of graphs achieving this equality was not given in \cite{gil-pons-2019}. We next settle this issue.

\begin{theorem}
\label{thm:dim=n-1}
A graph $G$ satisfies $\dimm(G) = n(G) - 1$ if and only if $G$ is a regular graph with $\diam(G)\le 2$.
\end{theorem}

\proof
It was observed in~\cite{gil-pons-2019} that $\dimm(K_n) = n-1$. Moreover, from~\cite[Example 3.3, Proposition 3.5]{gil-pons-2019} we also know that $\dimm(C_4) = 3$, $\dimm(C_5) = 4$, and $\dimm(C_n) = 3$ for $n\ge 6$, as well as that $\dimm(P_n)=1$. It follows that the theorem holds for complete graphs and graphs $G$ with $\Delta(G) \le 2$. In the rest of the proof we may thus assume that $\diam(G) \ge 2$ and $\Delta(G)\ge 3$.

Suppose first that $G$ is a graph which satisfies $\dimm(G) = n(G) - 1$. Assume that $G$ is not regular and select vertices $u$ and $v$ of $G$ such that $\deg_G(u) \ne \deg_G(v)$. We claim that $V(G)\setminus \{u, v\}$ is an outer multiset resolving set. For this sake we only need to verify that ${\rm m}_G(u|S) \ne {\rm m}_G(v|S)$. Let ${\rm m}_G(u|S) = \msl 1^{s_u}, \ldots\msr$ and ${\rm m}(v|S) = \msl 1^{s_v}, \ldots \msr$. Then, no matter whether $u$ and $v$ are adjacent or not, $\deg_G(u) \ne \deg_G(v)$ implies that $s_u \ne s_v$ and hence ${\rm m}_G(u|S) \ne {\rm m}_G(v|S)$. This shows that $G$ must be regular.

Suppose now that $G$ is an $r$-regular graph, $r\ge 3$, which satisfies $\dimm(G) = n(G) - 1$. We first claim that every vertex of $G$ is diametral. Assume on the contrary that there exists a non-diametral vertex $u$, and let $v$ be an arbitrary diametral vertex of $G$. Let $S = V(G) \setminus \{u,v\}$ and note that $\diam(G)\notin {\rm m}(u|S)$ while $\diam(G)\in {\rm m}(v|S)$. This means that $S$ is an outer multiset resolving set which in turn implies that $\dimm(G) < n(G) - 1$. This contradiction proves the claim that each vertex of $G$ is diametral.

Assume that $\diam(G) \ge 3$. Let $u$ be an arbitrary vertex of $G$ and consider a diametral path starting in $u$ and ending in $u'$. Let $w$ be the neighbor of $u$ on the diametral path.  Then $w\ne u'$ because $\diam(G) \ge 3$. Let $S = V(G)\setminus \{u,w\}$. Since $\dimm(G) = n(G) - 1$, the set $S$ is not an outer multiset resolving set. In particular, if ${\rm m}_G(u|S) = \msl \ldots, \diam(G)^{s_u}\msr$ and ${\rm m}_G(v|S) = \msl \ldots, \diam(G)^{s_w}\msr$, then $s_u = s_w > 0$. Consider now the set $T = V(G)\setminus \{u,u', w\}$. Note first that (having in mind that $\diam(G) \ge 3$ and that $G$ is $r$-regular) each of the vertices $u$ and $w$ has $r-1$ neighbors in $T$ while $u'$ has $r$ neighbors in $T$. Therefore, the pairs $u, u'$ and $w,u'$ are resolved by $T$. Moreover, since $u'\notin T$, there are $s_u-1$ vertices in $T$ at distance $\diam(G)$ from $u$, while there are $s_w$ vertices in $T$ that are at distance $\diam(G)$ from $w$. Since $s_u = s_w$ we have proved that $T$ also resolves the pair $u,w$. It follows that $T$ is an outer multiset resolving set, a contradiction to the assumption $\dimm(G) = n(G) - 1$. We came to this contradiction because we assumed that $\diam(G) \ge 3$. We conclude that $\diam(G) = 2$.

Conversely, suppose that $G$ is an $r$-regular graph, $r\ge 3$, with $\diam(G) = 2$. Let $S$ be an outer multiset basis of $G$. Assume that $|S| \le n(G) - 2$. Let $\overline{S} = V(G)\setminus S$. Since $\diam(G) = 2$, for each $x\in \overline{S}$ we have ${\rm m}_G(x|S) = \msl 1^{s_x}, 2^{s_x'}\msr$, where $s_x + s_x' = |S|$. As $S$ is an outer multiset resolving set, $s_x\ne s_y$ for each pair of vertices $x,y\in \overline{S}$, for otherwise $s_x = s_y$ would mean that $s'_x = s'_y$ and hence $x$ and $y$ have the same multiset representation. Select a vertex $u\in \overline{S}$ such that $s_u$ is smallest possible. Then, because $s_x\ne s_y$ for each pair of vertices $x,y\in \overline{S}$, and since $|\overline{S}| \ge 2$, we have $s_u < r$. Let $r = s_u + t$. Then $u$ has $t$ neighbors in $\overline{S}$. Because the vertices $x$ from $\overline{S}$ have pairwise different values $s_x$, and since $s_u$ is the smallest among them,  $\{s_x:\ x\in \overline{S}\} = \{s_u, s_u + 1, \ldots, s_u + t\} = \{s_u, s_u + 1, \ldots, r\}$. Hence, there exits a vertex $w\in \overline{S}$ with $s_w = r$. This means that $w$ has $r$ neighbors in $S$. But $w$ is adjacent also to $u\in \overline{S}$, hence $\deg_G(w) \ge r+1$. As this is not possible we conclude that $|S| = n(G) -1$, that is, $\dimm(G) = n(G) - 1$.
\qed

Examples of graphs from Theorem~\ref{thm:dim=n-1} are the Petersen graph, Hamming graphs (alias Cartesian products of complete graphs), and direct products of complete graphs. In addition, from Theorem~\ref{thm:dim=n-1} we immediately derive~\cite[Proposition 3.4]{gil-pons-2019} which asserts that if $k\ge 2$, then for the complete $k$-partite graph $K_{r, \ldots, r}$ we have $\dimm(K_{r, \ldots, r}) = kr-1$. On the other hand, we can easily get that if $k\ge 2$ and $2\le r_1 < r_2 < \cdots < r_k$, them $\dimm(K_{r_1, \ldots, r_k}) = r_1 + \cdots + r_k - k$.

\section{Graphs with outer multiset dimension $2$}
\label{sec:dim-2}

The problem of characterizing the graphs having its classical metric dimension equal to $2$ is one of the open problems in the area. See \cite{Behtoei} for an example with partial contributions in this direction. In this section we center our attention into those graphs $G$ with $\dimm(G) = 2$. We propose a polynomial algorithm for their recognition and describe their structure.

\begin{lemma}
\label{lem:mult-basis-2-dist-2}
If $G$ is a graph with $\dimm(G) = 2$ and $S=\{u,v\}$ is an outer multiset basis, then $d_G(u,v)\le 2$.
\end{lemma}

\proof
Suppose that $d_G(u,v)>2$. Let $P$ be a $u,v$-geodesic and let $u'\in N(u)\cap V(P)$ and $v'\in N(v)\cap V(P)$. Then $u'$ and $v'$ have the same multiset representation with respect to $S$, which is not possible.
\qed

Lemma~\ref{lem:mult-basis-2-dist-2} leads to Algorithm~\ref{alg:multiset-dim-2} which  decides in polynomial time whether $\dimm(G) = 2$ holds for a given graph $G$.

\begin{algorithm*}
\caption{Deciding whether a graph $G$ satisfies $\dimm(G) = 2$}
\label{alg:multiset-dim-2}
\begin{algorithmic}[1]
\Procedure{OUTER-MULTISET-DIMENSION-EQUAL-TWO}{$G$}
\If {$G$ is a path}
	\State \textbf{return} $\dimm(G) = 1$
\Else \State compute the distance matrix of $G$
\EndIf
\ForAll {$u \in V(G)$}
    \ForAll {$v \in V(G)$\;:\; $d_G(u,v)\le 2$}
        \If {multisets in $\{\msl d_G(u,x),d_G(v,x)\msr\,:\,x\ne u,v\}$ are different}
            \State $\dimm(G) = 2$
        \Else { }  \textbf{return} $\dimm(G) > 2$
        \EndIf
    \EndFor
\EndFor
\EndProcedure
\end{algorithmic}
\end{algorithm*}

\begin{theorem}
Deciding whether a graph $G$ of order $n$ satisfies $\dimm(G) = 2$ can be done in $\mathcal{O}(n^3)$ time.
\end{theorem}

\proof
By Lemma~\ref{lem:mult-basis-2-dist-2}, we only need to verify each pair of vertices at distance at most $2$ whether it forms an outer multiset basis. This checking is implemented in Algorithm~\ref{alg:multiset-dim-2} whose correctness it thus guaranteed by Lemma~\ref{lem:mult-basis-2-dist-2}.

The distance matrix of the graph $G$ of order $n$ and size $m$ can be computed in time $\mathcal{O}(nm)$. Moreover, since we have shown above that the maximum degree of $G$ is bounded by a (small) constant, $\mathcal{O}(nm) = \mathcal{O}(n^2)$. (We can pre-process $G$ by checking the degrees of its vertices.) The main loop (Step 6) is performed $n$ times, while for each vertex $u$, the inner loop (Step 7) in performed a constant number of times because the maximum degree of $G$ and hence also the square of the maximum degree are constant. Checking whether multisets in Step 8 are different can be done in $\mathcal{O}(n^2)$ time, hence the total complexity is $\mathcal{O}(n^3)$.
\qed

Despite the fact that deciding whether the outer multiset dimension of a graph equals $2$ is polynomial, it is of interest to have more insight into the structure of such graphs. To do so, for a graph $G$, $X\subseteq V(G)$, and $k\ge 0$, we define
$$L_k(X) = \{u\in V(G):\ \min_{x\in X} d_G(u,x) = k\}.$$
Note that $L_0(X) = X$, and that the sets $L_k(X)$, $k\ge 0$, partition $V(G)$.

\begin{lemma}
\label{lem:layers}
Let $G$ be a graph with $\dimm(G) = 2$ and let $S$ be an outer multiset basis of $G$. Then for every $k\ge 1$ we have $|L_k(S)| \le 3$. Moreover, if the vertices of $S$ are adjacent, then $|L_k(S)| \le 2$ and $|L_{k+1}(S)| \le |L_k(S)|$.
\end{lemma}

\proof
Let $S = \{u, v\}$. Let $k\ge 1$ and let $x\in L_k(S)$. By definition of $L_k(S)$, we may without loss of generality assume that $d_G(x,u) = k$.  Since $d_G(u,v)\le 2$ by Lemma~\ref{lem:mult-basis-2-dist-2}, we have $d_G(v,x) \in \{k, k+1, k+2\}$. It follows that $m(x|S)$ is one of $\msl k, k\msr$, $\msl k, k+1\msr$, and $\msl k, k+2\msr$. As $S$ is an outer multiset  basis this in turn implies that $|L_k(S)| \le 3$.

Assume in the rest that $d_G(u,v) = 1$. Then $d_G(v,x) \in \{k, k+1\}$ and we can conclude similarly as above that $|L_k(S)| \le 2$. Suppose now that $|L_k(S)| = 1$ and $|L_{k+1}(S)| = 2$ for some $k\ge 1$. Let $L_{k+1}(S) = \{x, y\}$ and $L_{k}(S) = \{z\}$. Then, by the definition of the sets $L_i(S)$, we infer that $z$ is adjacent to both $x$ and $y$. But this means that $x$ and $y$ have the same mutiset representation with respect to $S$. This contradiction proves that $|L_{k+1}(S)| \le |L_k(S)|$.
\qed

Based on Lemma~\ref{lem:layers}, we next characterize the graphs $G$ with $\dimm(G) = 2$ having an outer multiset basis formed by two adjacent vertices. To this end, let $S = \{u,v\}$ be an outer multiset basis of $G$ with $uv\in E(G)$. By Lemma~\ref{lem:layers}, $|L_k(S)| \le 2$ for each $k$. Suppose now that for some $k$ we have $L_k(S) = \{x_k, y_k\}$ and $L_{k+1}(S) = \{x_{k+1}, y_{k+1}\}$. Assume without loss of generality that $m(x_k|S) = \msl k, k\msr$, $m(y_{k}|S) = \msl k, k+1\msr$, $m(x_{k+1}|S) = \msl k+1, k+1\msr$, and $m(y_{k+1}|S) = \msl k+1, k+2\msr$. Then $x_k$ is not adjacent to $y_{k+1}$, but must be adjacent to $x_{k+1}$. In addition, $y_k$ must be adjacent to $y_{k+1}$, and may be adjacent to $x_{k+1}$. Finally each of the edges $x_ky_{k}$ and $x_{k+1}y_{k+1}$ may be present.

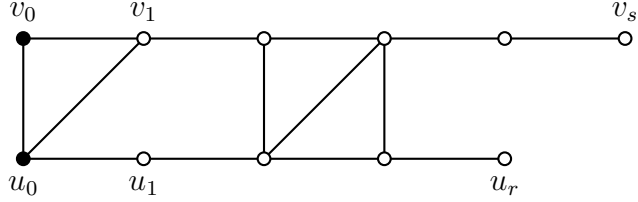
\begin{figure}[h!]
\begin{center}
\begin{tikzpicture}[scale=0.8,style=thick]
\tikzstyle{every node}=[draw=none,fill=none]
\def\vr{3pt} 

\begin{scope}[yshift = 0cm, xshift = 0cm]
\path (0,0) coordinate (v1);
\path (2,0) coordinate (v2);
\path (4,0) coordinate (v3);
\path (6,0) coordinate (v4);
\path (8,0) coordinate (v5);
\path (0,2) coordinate (v6);
\path (2,2) coordinate (v7);
\path (4,2) coordinate (v8);
\path (6,2) coordinate (v9);
\path (8,2) coordinate (v10);
\path (10,2) coordinate (v11);
\draw (v5) -- (v1) -- (v6) -- (v11);
\draw (v1) -- (v7);
\draw (v8)-- (v3) -- (v9) -- (v4);
\draw (v1)  [fill=black] circle (\vr);
\draw (v2)  [fill=white] circle (\vr);
\draw (v3)  [fill=white] circle (\vr);
\draw (v4)  [fill=white] circle (\vr);
\draw (v5)  [fill=white] circle (\vr);
\draw (v6)  [fill=black] circle (\vr);
\draw (v7)  [fill=white] circle (\vr);
\draw (v8)  [fill=white] circle (\vr);
\draw (v9)  [fill=white] circle (\vr);
\draw (v10)  [fill=white] circle (\vr);
\draw (v11)  [fill=white] circle (\vr);
\draw[below] (v1)++(0.0,-0.1) node {$u_0$};
\draw[below] (v2)++(0.0,-0.1) node {$u_1$};
\draw[below] (v5)++(0.0,-0.1) node {$u_r$};
\draw[above] (v6)++(0.0,0.1) node {$v_0$};
\draw[above] (v7)++(0.0,0.1) node {$v_1$};
\draw[above] (v11)++(0.0,0.1) node {$v_s$};
\end{scope}

\end{tikzpicture}
\end{center}
\vspace*{-0.5cm}
\caption{A fairly representative example of a graph from the family $\mathcal{F}$.}
\label{fig:dim-2-adjacent}
\end{figure}

The above description together with the fact of Lemma~\ref{lem:layers} that for each $k\ge 1$ we have $|L_{k+1}(S)| \le |L_k(S)|$, lead to the  family $\mathcal{F}$ of graphs $G$ constructed in the following way. The vertex set of $G\in\mathcal{F}$ is $V(G)=\{u_0,\dots,u_r\}\cup \{v_0,\dots,v_s\}$ for some $r\ge 0$ and $s\ge 1$, and the edges of $G$ are given as follows.
\begin{itemize}
  \item $u_0v_0,u_0v_1\in E(G)$.
  \item For every $i\in [r]$ and every $j\in [s]$, $u_{i-1}u_{i}\in E(G)$ and $v_{j-1}v_{j}\in E(G)$.
  \item For every $i\in [\min\{r,s\}]$, the edge $u_iv_i$ might exist or not in $G$.
  \item For every $j\in [\min\{r,s-1\}]$, the edge $u_iv_{i+1}$ might exist or not in $G$.
\end{itemize}
Note that for instance, if $r=0$ and $s=1$ in the construction, then $G$ is precisely the complete graph $K_3$. If $r=0$ and $s=2$, then $G$ is the {\em paw}, that is, the graph obtained from $K_3$ by attaching a pendant edge to one of its vertices; and if $r=s=1$, then $G$ is either $K_4$ minus an edge, or the paw. For a representative example of a graph from the family $\mathcal{F}$ see Fig.~\ref{fig:dim-2-adjacent}. The outer multiset basis is indicated bold.

The above discussion gives us the following result.

\begin{theorem}
\label{th:multiset-dim-2-adj}
A graph $G$ has an outer multiset basis formed by two adjacent vertices if and only if $G\in \mathcal{F}$.
\end{theorem}

\begin{figure}[ht!]
\begin{center}
\begin{tikzpicture}[scale=0.9,style=thick]
\tikzstyle{every node}=[draw=none,fill=none]
\def\vr{3pt} 

\begin{scope}[yshift = 0cm, xshift = 0cm]
\path (2,0) coordinate (v1);
\path (4,0) coordinate (v2);
\path (6,0) coordinate (v3);
\path (0,2) coordinate (v4);
\path (2,2) coordinate (v5);
\path (4,2) coordinate (v6);
\path (0,4) coordinate (v8);
\path (2,4) coordinate (v9);
\path (4,4) coordinate (v10);
\path (6,4) coordinate (v11);
\path (8,4) coordinate (v12);
\path (8,2) coordinate (v13);
\path (8,0) coordinate (v14);
\path (10,2) coordinate (v15);
\draw (v1)-- (v2) -- (v3);
\draw (v2) -- (v6) -- (v1) -- (v4) -- (v5) -- (v6) -- (v1) -- (v2);
\draw (v6) -- (v10) -- (v5) -- (v6) -- (v4) -- (v5) -- (v10) -- (v9) -- (v8);
\draw (v1)-- (v6) -- (v11);
\draw (v11) -- (v10);
\draw (v4) -- (v9) -- (v5);
\draw (v14) -- (v3) -- (v13);
\draw (v11) -- (v12);
\draw (v12) -- (v13) -- (v15);
\draw (v1)  [fill=white] circle (\vr);
\draw (v2)  [fill=white] circle (\vr);
\draw (v3)  [fill=white] circle (\vr);
\draw (v4)  [fill=black] circle (\vr);
\draw (v5)  [fill=white] circle (\vr);
\draw (v6)  [fill=white] circle (\vr);
\draw (v8)  [fill=black] circle (\vr);
\draw (v9)  [fill=white] circle (\vr);
\draw (v10)  [fill=white] circle (\vr);
\draw (v11)  [fill=white] circle (\vr);
\draw (v12)  [fill=white] circle (\vr);
\draw (v13)  [fill=white] circle (\vr);
\draw (v14)  [fill=white] circle (\vr);
\draw (v15)  [fill=white] circle (\vr);
{\small
\draw[above] (v9)++(0.0,0.1) node {$\msl 1,1\msr$};
\draw[above] (v10)++(0.0,0.1) node {$\msl 2,2\msr$};
\draw[above] (v11)++(0.0,0.1) node {$\msl 3,3\msr$};
\draw[above] (v12)++(0.0,0.1) node {$\msl 4,4\msr$};
\draw[below] (v5)++(0.0,-0.1) node {$\msl 1,2\msr$};
\draw[right] (v6)++(+0.0,-0.2) node {$\msl 2,3\msr$};
\draw[left] (v13)++(-0.1,0.0) node {$\msl 4,5\msr$};
\draw[above] (v15)++(0.0,0.1) node {$\msl 5,6\msr$};
\draw[below] (v1)++(0.0,-0.1) node {$\msl 1,3\msr$};
\draw[below] (v2)++(0.0,-0.1) node {$\msl 2,4\msr$};
\draw[below] (v3)++(0.0,-0.1) node {$\msl 3,5\msr$};
\draw[below] (v14)++(0.0,-0.1) node {$\msl 4,6\msr$};
}
\end{scope}

\end{tikzpicture}
\end{center}
\vspace*{-5mm}
\caption{A fairly representative example of a graph $G$ with $\dimm(G) = 2$ and an outer multiset basis formed by two non adjacent vertices.}
\label{fig:dim-2-not-adjacent}
\end{figure}
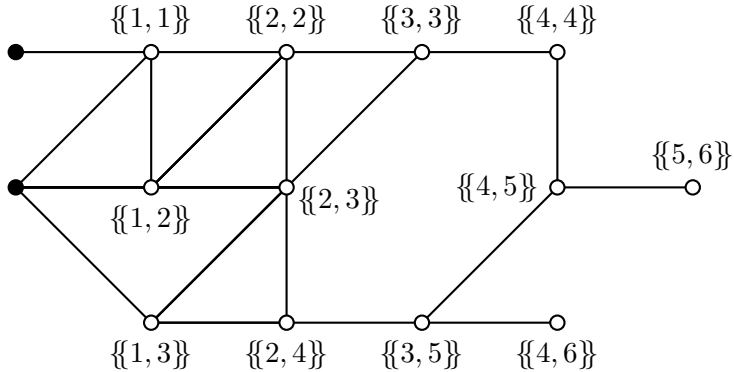

In view of Lemma~\ref{lem:mult-basis-2-dist-2} and Theorem~\ref{th:multiset-dim-2-adj}, in order to complete the characterization of the graphs $G$ with $\dimm(G) = 2$, it remains to describe how such graphs look like for the case in which all their outer multiset bases are formed by two non adjacent vertices. This can be done in a similar way as we did it when there exists an outer multiset basis consisting of two adjacent vertices, however a formal description is more technical and hence not given in detail. The reason for this is the fact that in this subcase it is not satisfied the property $|L_{k+1}(S)| \le |L_k(S)|$ which holds when there exists an outer mutiset basis with two adjacent vertices. Nevertheless, we still have that $|L_k(S)|\le 3$. Thus, instead of presenting a lengthy description, we only show a typical example in Fig.~\ref{fig:dim-2-not-adjacent} with its outer multiset basis again in bold.

\section{Lexicographic products and multiset distance irregular graphs}
\label{sec:lexico}

In this section we consider the outer multiset dimension of lexicographic products of graphs. Recall that the {\em lexicographic product} $G\circ H$ of graphs $G$ and $H$ has the vertex set $V(G)\times V(H)$ and edges $(g,h)(g',h')$, where either $g=g'$ and $hh'\in E(H)$, or $gg'\in E(G)$. If $g\in V(G)$, then the set of vertices $\{(g,h):\ h\in V(H)\}$ induces a subgraph of $G\circ H$ isomorphic to $H$ which is called an {\em $H$-layer} and denoted by $^{g}H$.

The metric dimension of lexicographic products has been independently investigated in papers~\cite{jannesari-2012, saputro-2013}. The main tool used in these articles is that one of transforming the metric dimension in the lexicographic product $G\circ H$ to the so-called adjacency dimension of $H$. Indeed, the adjacency dimension was explicitly introduced for the first time in~\cite{jannesari-2012} (and implicitly in~\cite{saputro-2013}), for more information on it see~\cite{dorota-2022+}.

We say that a graph $G$ is {\em multiset distance irregular} if for every two vertices $u,v\in V(G)$, the multisets ${\rm m}_G(u|V(G))$ and ${\rm m}_G(v|V(G))$ are different.

\begin{theorem}
\label{thm:lex}
If $G$ is a graph with $n(G)\ge 2$ and $H\in \{K_k, \overline{K_k}\}$, $k\ge 2$, then
$$\dimm(G\circ H) \ge n(G)(k-1).$$
Moreover, equality holds if and only if $G$ is multiset distance irregular.
\end{theorem}

\proof
Let $V(G) = \{g_1, \ldots, g_{n(G)}\}$. For $g_i$, $i\in [n(G)]$, let $d_i = \ecc_G(g_i)$ be the {\em eccentricity} of $g_i$ which is the largest distance between $g_i$ and any other vertex of $G$. Then write
\begin{equation}
\label{eq:in-G}
{\rm m}_G(g_i|V(G)) = \msl 0, 1^{i_1}, \ldots, d_i^{i_d} \msr.
\end{equation}
Let $(g_i,h)$ and $(g_i,h')$ be different vertices from the $^{g_i}H$-layer. Since $H\in \{K_k, \overline{K_k}\}$, the vertices $(g_i,h)$ and $(g_i,h')$ are twins. By~\cite[Proposition 3.7]{gil-pons-2019} we know that every outer multiset resolving set of $G\circ H$ contains at least one of $(g_i,h)$ and $(g_i,h')$. Inductively, every outer multiset resolving set of $G\circ H$ contains at least $k-1$ vertices from $^{g_i}H$. As this $H$-layer was arbitrary, the inequality follows.

Suppose now that the equality holds and let $S$ be an outer multiset basis of $G\circ H$. By the above argument, each layer $^{g_i}H$ contains exactly one vertex which does not belong to $S$. We may without loss of generality assume that this vertex is $(g_i,h)$, where $h$ is some fixed vertex of $H$. In view of~\eqref{eq:in-G} we then have
\begin{equation}
\label{eq:complete}
{\rm m}_{G\circ H}((g_i,h)|S) = \msl 1^{i_1(k-1)}, \ldots, d_i^{i_d(k-1)}, 1^{(k-1)} \msr
\end{equation}
when $H = K_k$ and
\begin{equation}
\label{eq:non-complete}
{\rm m}_{G\circ H}((g_i,h)|S) = \msl 1^{i_1(k-1)}, \ldots, d_i^{i_d(k-1)}, 2^{(k-1)} \msr
\end{equation}
when $H = \overline{K_k}$. The above respective terms $1^{(k-1)}$ and $2^{(k-1)}$ reflect the distances between $(g_i,h)$ and the other vertices from $^{g_i}H$ (which are all in $S$). Since $S$ is an outer multiset basis of $G\circ H$, we have ${\rm m}_{G\circ H}((g_i,h)|S) \ne {\rm m}_{G\circ H}((g_j,h)|S)$ for each $i,j\in [n(G)]$, $i\ne j$.
But then no matter whether we consider~\eqref{eq:complete} or~\eqref{eq:non-complete}, we deduce that ${\rm m}_{G}(g_i|V(G)) \ne {\rm m}_{G}(g_j|V(G))$ holds for each $i,j\in [n(G)]$. Hence $G$ is a multiset distance irregular.

Conversely, let $G$ be a multiset distance irregular graph. Then, by definition, ${\rm m}_{G}(g_i|V(G)) \ne {\rm m}_{G}(g_j|V(G))$ holds for each $i,j\in [n(G)]$. Let $S\subseteq V(G\circ H)$ be a set that contains exactly $k-1$ vertices of each $H$-layer. But then in view of \eqref{eq:complete} or~\eqref{eq:non-complete}, ${\rm m}_{G\circ H}((g_i,h)|S) \ne {\rm m}_{G\circ H}((g_j,h)|S)$ for each $i,j\in [n(G)]$, $i\ne j$. This means that $S$ is a outer multiset resolving set. Since $|S| = n(G)(k-1)$ we are done.
\qed

We believe that multiset distance irregular graphs are of independent interest. First of all, they are closely related to transmission irregular graphs which are defined as follows. The {\em transmission} $\Tr_G(v)$ of a vertex $v$ of a graph $G$ is the sum of distances between $v$ and all the other vertices of the graph. $G$ is {\em transmission irregular} if the vertices of $G$ have pairwise different transmissions. (Transmission irregular graphs are also known as status injective graphs, cf.~\cite{QZ2020}.) This concept was first studied in~\cite{AK2018} with respect to the Wiener dimension because transmission irregular graphs are the graphs with a largest possible Wiener dimension. For a selection of different appealing constructions of transmission irregular graphs and results on these graphs see~\cite{al-yakoob-2020, dobrynin-2019a, xu-2021}.

\begin{figure}[ht!]
\begin{center}
\begin{tikzpicture}[scale=0.9,style=thick]
\tikzstyle{every node}=[draw=none,fill=none]
\def\vr{3pt} 

\begin{scope}[yshift = 0cm, xshift = 0cm]
\path (0,0) coordinate (v1);
\path (3,0) coordinate (v2);
\path (6,0) coordinate (v3);
\path (9,0) coordinate (v4);
\path (12,0) coordinate (v5);
\path (3,3) coordinate (v6);
\path (6,3) coordinate (v7);
\path (9,3) coordinate (v8);
\draw (v1) -- (v2) -- (v3) -- (v4) -- (v5);
\draw (v2) -- (v6) -- (v7) -- (v8);
\draw (v3) -- (v7) -- (v4);
\draw (v1)  [fill=white] circle (\vr);
\draw (v2)  [fill=white] circle (\vr);
\draw (v3)  [fill=white] circle (\vr);
\draw (v4)  [fill=white] circle (\vr);
\draw (v5)  [fill=white] circle (\vr);
\draw (v6)  [fill=white] circle (\vr);
\draw (v7)  [fill=white] circle (\vr);
\draw (v8)  [fill=white] circle (\vr);
\draw[below] (v1)++(0.0,-0.1) node {$\msl 1^1, 2^2, 3^2, 4^2 \msr$};
\draw[below] (v2)++(0.0,-0.1) node {$\msl 1^3, 2^2, 3^2 \msr$};
\draw[below] (v3)++(0.0,-0.1) node {$\msl 1^3, 2^4 \msr$};
\draw[below] (v4)++(0.0,-0.1) node {$\msl 1^3, 2^3, 3^1 \msr$};
\draw[below] (v5)++(0.0,-0.1) node {$\msl 1^1, 2^2, 3^3, 4^1 \msr$};
\draw[above] (v6)++(0.0,0.1) node {$\msl 1^2, 2^4, 3^1 \msr$};
\draw[above] (v7)++(0.0,0.1) node {$\msl 1^4, 2^2, 3^1 \msr$};
\draw[above] (v8)++(0.0,0.1) node {$\msl 1^1, 2^3, 3^2, 4^1 \msr$};
\draw[above] (v1)++(0.0,0.1) node {$19$};
\draw[above] (v2)++(0.3,0.1) node {$13$};
\draw[above] (v3)++(-0.3,0.1) node {$11$};
\draw[above] (v4)++(0.2,0.1) node {$12$};
\draw[above] (v5)++(0.0,0.1) node {$18$};
\draw[below] (v6)++(0.3,-0.1) node {$13$};
\draw[below] (v7)++(-0.3,-0.1) node {$11 $};
\draw[below] (v8)++(0.0,-0.1) node {$17$};
\end{scope}

\end{tikzpicture}
\end{center}
\vspace{-5mm}
\caption{A graph which is multiset distance irregular but not transmission irregular.}
\label{fig:example}
\end{figure}
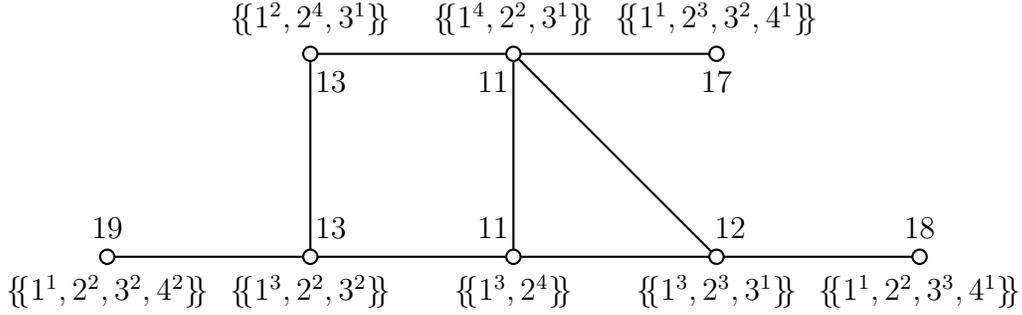

Clearly, if ${\rm m}_{G}(u|V(G)) = {\rm m}_{G}(v|V(G))$, then $\Tr_G(u) = \Tr_G(v)$. Hence transmission irregular graphs form a subset of multiset distance irregular graphs. The inclusion is strict as demonstrated by the graph $X$ from Fig.~\ref{fig:example} which is multiset distance irregular but not transmission irregular. Next to each vertex $u$ the multiset ${\rm m}_G(u|(V(X)\setminus\{u\})$ is written as well as $\Tr_X(u)$. The multisets ${\rm m}_G(u|V(X))$ are indeed pairwise different, but there are two pairs of vertices with the same transmission.

\section{Grid graphs}
\label{sec:cartesian}

Let $V(P_n)=\{0,1,\dots,n-1\}$ where $i$ is adjacent to $j$ if and only if $|i-j|=1$. The grid graph $P_{s}\,\Box\, P_{t}$ is the Cartesian product of two paths $P_s$ and $P_t$, that is $V(P_{s}\,\Box\, P_{t})=\{(i,j)\,:\,0\le i\le s-1\mbox{ and } 0\le j\le t-1\}$ and $(i,j)(k,\ell)\in E(P_{s}\,\Box\, P_{t})$ when $|i-j|+|k-\ell|=1$.

The metric dimension of grids were reported in several papers, among the earliest ones are~\cite{khuller-1996, melter-1984}. Here we add the following result for the outer metric dimension of grids.

\begin{theorem}
\label{thm:grid}
For all $s\ge t\ge 2$,
$$\dimm(P_s\Box P_t)=3.$$
\end{theorem}

\proof
In order to simplify the notation, let $G=P_s\Box P_t$. It clearly cannot be $\dimm(G)=1$ since $G$ is not a path. Also, notice that if $s=t=2$, then by Theorem \ref{thm:dim=n-1} we have that $\dimm(G)=3.$ Accordingly, from now on we may w.l.g. assume that $s\ge 3$. If $s=3$, then $t\in \{2,3\}$. In such a case it can be easily checked that the set $\{(0,0), (2,0), (2,1)\}$ is an outer multiset basis. Hence, from now we assume that $s\ge 4$.

Now, if $\dimm(G)=2$ and $S=\{(i,j),(i',j')\}$ is an outer multiset basis, then by Lemma~\ref{lem:mult-basis-2-dist-2}, we have that $d_{G}((i,j),(i',j'))\le 2$. If $i=i'$, we readily observe that the two vertices $(i+1,j),(i+1,j')$ (or the two vertices $(i-1,j),(i-1,j')$ if $i=s$) have the same multiset representation with respect to $S$, which is not possible. By symmetry, a similar conclusion is deduced if $j=j'$. It remains the case $d_{G}((i,j),(i',j'))=2$, $i\ne i'$ and $j\ne j'$. Thus, it must happen that (w.l.g.) $S=\{(i,j),(i+1,j+1)\}$. Hence, $(i+1,j)$ and $(i,j+1)$ have the same multiset representation with respect to $S$, and this is also not possible. As a consequence of the two previous contradictions, we obtain that $\dimm(G)\ge 3$.

In order to show that $\dimm(G)\le 3$, we claim that $R=\{(0,0),(1,0),(s-1,0)\}$ is an outer multiset resolving set for $G$. Let $R'=\{(0,0),(s-1,0)\}$ and notice that $R'$ is a metric basis for $G$, cf.~\cite{khuller-1996}. Let $(i,j),(i',j')$ be any two vertices in $V(G)\setminus R$. We observe that the multiset representations of $(i,j)$ and $(i',j')$ with respect to $R'$ are:
\begin{equation*}
{\rm m}_G((i,j)|R')=\msl i+j,s-1-i+j \msr
\end{equation*}
\begin{equation*}
{\rm m}_G((i',j')|R')=\msl i'+j',s-1-i'+j' \msr.
\end{equation*}
Since $R'$ is a metric basis, we have one of the following three situations.

\medskip
\noindent
Case 1: $d_G((i,j),(0,0))=d_G((i',j'),(0,0))=\alpha$. \\
In such situation, it must happen $\beta = d_G((i,j),(s-1,0))\ne d_G((i',j'),(s-1,0)) = \beta'$ since $R'$ is a metric basis. On the other hand, since $G$ is bipartite,
$$d_G((i,j),(1,0)) \in \{d_G((i,j),(0,0))-1, d_G((i,j),(0,0))+1\},$$
and
$$d_G((i',j'),(1,0)) \in \{d_G((i',j'),(0,0))-1, d_G((i',j'),(0,0))+1\}.$$
If $d_G((i,j),(1,0))=d_G((i,j),(0,0))+1=\alpha+1$ and $d_G((i',j'),(1,0))=d_G((i',j'),(0,0))+1=\alpha+1$, then we note that
$${\rm m}_G((i,j)|R)=(\alpha,\beta,\alpha+1) \ne (\alpha,\beta', \alpha+1)={\rm m}_G((i',j')|R)$$
since $\beta\ne \beta'$. Similarly, if $d_G((i,j),(1,0))=d_G((i,j),(0,0))-1$ and $d_G((i',j'),(1,0))=d_G((i',j'),(0,0))-1$, then we note that ${\rm m}_G((i,j)|R)\ne {\rm m}_G((i',j')|R)$.

We now consider the subcase when $d_G((i,j),(1,0))= \alpha+1$ and $d_G((i',j'),(1,0)) = \alpha-1$. From the first equality, we deduce that $i=0$ must occur. Thus, from the case assumption we obtain that $j=i+j=d_G((i,j),(0,0))=d_G((i',j'),(0,0))=i'+j'$. Notice that the multiset representations of $(i,j)$ and $(i',j')$ with respect to $R$ are:
$${\rm m}_G((i,j)|R)=\msl i+j,s-1-i+j,i+j+1 \msr=\msl j,s-1+j,j+1 \msr,$$
and
$${\rm m}_G((i',j')|R)=\msl i'+j',s-1-i'+j',i'+j'-1 \msr.$$
Recall that $j=i'+j'$ and that $s-1+j\ne s-1-i'+j'$. If ${\rm m}_G((i,j)|R)={\rm m}_G((i',j')|R)$, then it must happen that $s-1+j=i'+j'-1$ and that $j+1=s-1-i'+j'$. However, since $j=i'+j'$, from $s-1+j=i'+j'-1$ we deduce that $s=0$, which is not possible. Thus, ${\rm m}_G((i,j)|R)\ne{\rm m}_G((i',j')|R)$.

To complete this case, we may assume $d_G((i,j),(1,0)) = \alpha-1$ and $d_G((i',j'),(1,0)) = \alpha +1$. We again obtain a similar conclusion, by using a procedure analogous to the one above, but taking into account that now it must first occur $i'=0$ instead of $i=0$.

\medskip
\noindent
Case 2: $d_G((i,j),(s-1,0))=d_G((i',j'),(s-1,0))$. \\In such situation, it must happen $d_G((i,j),(0,0))\ne d_G((i',j'),(0,0))$ since $R'$ is a metric basis. If ${\rm m}_G((i,j)|R)={\rm m}_G((i',j')|R)$, then by using some similar arguments as in Case 1, it must happen that the multisets of $(i,j)$ and $(i',j')$ with respect to $R$ are either
$${\rm m}_G((i,j)|R)=\msl i+j,s-1-i+j,i+j+1 \msr,$$
and
$${\rm m}_G((i',j')|R)=\msl i'+j',s-1-i'+j',i'+j'-1 \msr;$$
or
$${\rm m}_G((i,j)|R)=\msl i+j,s-1-i+j,i+j-1 \msr,$$
and
$${\rm m}_G((i',j')|R)=\msl i'+j',s-1-i'+j',i'+j'+1 \msr.$$
Moreover, $i+j=i'+j'-1$ and $i'+j'=i+j+1$ must happen in the first situation, as well as $i+j=i'+j'+1$ and $i'+j'=i+j-1$ in the second one.

In addition, in the first possibility we have that $i=0$ must happen, while in the second one $i'=0$ must happen instead. Since $s-1-i+j=d_G((i,j),(s-1,0))=d_G((i',j'),(s-1,0))=s-1-i'+j'$, we obtain that either $j=j'-i'$ (when $i=0$), or $j'=j-i$ (when $i'=0$). By using these equalities in $i+j=i'+j'-1$ and $i+j=i'+j'+1$, we deduce that $2i'=-1$ and that $2i'=1$, respectively, which are both contradictions. Consequently, we again obtain ${\rm m}_G((i,j)|R)\ne {\rm m}_G((i',j')|R)$.

\medskip
\noindent
Case 3: $d_G((i,j),(0,0))\ne d_G((i',j'),(0,0))$ and $d_G((i,j),(s-1,0))\ne d_G((i',j'),(s-1,0))$.\\
Clearly, if $d_G((i,j),(0,0))\ne d_G((i',j'),(s-1,0))$ and $d_G((i,j),(s-1,0))\ne d_G((i',j'),(0,0))$, then we have that ${\rm m}_G((i,j)|R)\ne {\rm m}_G((i',j')|R)$, independently on which the distances $d_G((i,j),(1,0))$ and $d_G((i',j'),(1,0))$ are.

Now, if $$i+j=d_G((i,j),(0,0))=d_G((i',j'),(s-1,0))=s-1-i'+j'$$
and
$$s-1-i+j=d_G((i,j),(s-1,0))=d_G((i',j'),(0,0))=i'+j',$$
then we deduce that $j=j'$. Thus, since $s\ge 4$ it must happen that $d_G((i,j),(1,0))\ne d_G((i',j'),(1,0))$, for otherwise we get either $i=i'$ or (w.l.g.) $i=0$ and $i'=2$ (which are contradictions). This leads to ${\rm m}_G((i,j)|R)\ne {\rm m}_G((i',j')|R)$.

To conclude the proof, we need to consider the case when
$$i+j=d_G((i,j),(0,0))=d_G((i',j'),(s-1,0))=s-1-i'+j'$$
and
$$s-1-i+j=d_G((i,j),(s-1,0))\ne d_G((i',j'),(0,0))=i'+j'.$$
If ${\rm m}_G((i,j)|R)={\rm m}_G((i',j')|R)$, then it must happen $d_G((i,j),(s-1,0))=d_G((i',j'),(1,0))$ and that $d_G((i,j),(1,0))=d_G((i',j'),(0,0))$. Moreover, since
$$d_G((i,j),(1,0))= d_G((i,j),(0,0))+1\mbox{ or }d_G((i,j),(1,0))= d_G((i,j),(0,0))-1,$$ and
$$d_G((i',j'),(1,0))=d_G((i',j'),(0,0))+1\mbox{ or }d_G((i',j'),(1,0))=d_G((i',j'),(0,0))-1,$$
it must be satisfied that $|d_G((i,j),(s-1,0))- d_G((i',j'),(0,0))|=1$. However, by using the fact that $i+j=d_G((i,j),(0,0))=d_G((i',j'),(s-1,0))=s-1-i'+j'$ in the difference $d_G((i,j),(s-1,0))-d_G((i',j'),(0,0))=(s-1-i+j)-(i'+j')$, we obtain that $d_G((i,j),(s-1,0))-d_G((i',j'),(0,0))=2(j-j')$, which is an even number, a contradiction. Therefore, ${\rm m}_G((i,j)|R)={\rm m}_G((i',j')|R)$. A similar conclusion is obtained if $i+j=d_G((i,j),(0,0))\ne d_G((i',j'),(s-1,0))=s-1-i'+j'$ and $s-1-i+j=d_G((i,j),(s-1,0))=d_G((i',j'),(0,0))=i'+j'$.

As a consequence of the arguments above, we conclude that $R$ is an outer multiset resolving set, which completes the proof.
\qed

\section{Concluding remarks}

The discussion at the end of Section~\ref{sec:minus-one} indicates that it would be of interest to investigate the outer multiset dimension of non-regular graphs of diameter $2$. In Section~\ref{sec:lexico}, a lower bound is proved for the lexicographic products in which the second factor is complete or edgeless. Studying lexicographic products in general is an open area. Another challenge is to extend Theorem~\ref{thm:grid} to multidimensional grids, that is, to Cartesian products of several paths. In this direction, hypercubes (that is, Cartesian products of paths of order $2$) deserve a special attention. An additional  class for which it would be interesting to determine the outer multiset dimension is the class of torus graphs $C_s\cp C_t$, $s,t\ge 3$.

\newpage
\section*{Acknowledgements}

Sandi Klav\v{z}ar acknowledges the financial support from the Slovenian Research Agency through research core funding No.\ P1-0297 and projects J1-2452 and N1-0285. Dorota Kuziak and Ismael G. Yero have been partially supported by the Spanish Ministry of Science and Innovation through the grant PID2019-105824GB-I00. Moreover, this work was initiated while the first author Sandi Klav\v{z}ar was visiting the University of Cadiz with the support of the grant PID2019-105824GB-I00. Dorota Kuziak has also  been partially supported by ``Plan Propio de Investigaci\'{o}n-UCA" ref.\ no.\ EST2022-075.


\begin{thebibliography}{99}

\bibitem{AK2018}
  Y.~Alizadeh, S.~Klav\v{z}ar,
  On graphs whose Wiener complexity equals their order and on Wiener index of asymmetric graphs,
  Appl.\ Math.\ Comput.\ 328 (2018) 113--118.

\bibitem{al-yakoob-2020}
  S.~Al-Yakoob, D.~Stevanovi\'{c},
  On transmission irregular starlike trees,
  Appl.\ Math.\ Comput.\ 380 (2020) 125257.

\bibitem{Behtoei}
  A.~Behtoei, A.~Davoodi, M.~Jannesari, B.~Omoomi,
  A characterization of some graphs with metric dimension two,
  Discrete Math. Algorithms Appl. 9 (2017) 1750027.

\bibitem{bong-2021}
  N.~H.~Bong, Y.~Lin,
  Some properties of the multiset dimension of graphs,
  Electron.\ J.\ Graph Theory Appl.\ 9 (2021) 215--221.

\bibitem{dobrynin-2019a}
  A.~A.~Dobrynin,
  Infinite family of transmission irregular trees of even order,
  Discrete Math.\ 342 (2019) 74--77.

\bibitem{gil-pons-2019}
  R.~Gil-Pons, Y.~Ram\'{\i}rez-Cruz, R.~Trujillo-Rasua, I.~G.~Yero,
  Distance-based vertex identification in graphs: the outer multiset dimension,
  Appl.\ Math.\ Comput.\ 363 (2019) 124612.

\bibitem{Harary}
  F. Harary, R.~A.~Melter,
  On the metric dimension of a graph,
  Ars Combin.\ 2 (1976) 191--195.

\bibitem{jannesari-2012}
  M.~Jannesari, B.~Omoomi,
  The metric dimension of the lexicographic product of graphs,
  Discrete Math.\ 312 (2012) 3349--3356.

\bibitem{khuller-1996}
  S.~Khuller B.~Raghavachari, A.~Rosenfeld,
  Landmarks in graphs,
  Discrete Appl.\ Math.\ 70 (1996) 217--229.

\bibitem{dorota-2022+}
  D.~Kuziak, I.~G.~Yero,
  Metric dimension related parameters in graphs: A survey on combinatorial, computational and applied results,
  arXiv:2107.04877 [math.CO] (10 Jul 2021).

\bibitem{melter-1984}
  R.~A.~Melter, I.~Tomescu,
  Metric bases in digital geometry,
  Comput.\ Vision Graphics Image Process.\ 25 (1984) 113--121.

\bibitem{QZ2020}
  P.~Qiao, X.~Zhan,
  Pairs of a tree and a nontree graph with the same status sequence,
  Discrete Math.\  343 (2020) 111662.

\bibitem{saputro-2013}
  S.~W.~Saputro, R.~Simanjuntak, S.~Uttunggadewa, H.~Assiyatun, E.~T.~Baskoro, A.~N.~M.~Salman, M.~Ba\v{c}a,
  The metric dimension of the lexicographic product of graphs,
  Discrete Math.\ 313 (2013) 1045--1051.

\bibitem{rino-2017}
  R.~Simanjuntak, T.~Vetr\'{\i}k, P.~Bintang Mulia,
  The multiset dimension of graphs,
  arXiv:1711.00225 [math.CO] (1 Nov 2017).

\bibitem{Slater}
  P.~J.~Slater,
  Leaves of trees,
  Cong.\ Numer.\ 14 (1975) 549--559.

\bibitem{till-2022+}
  R.~C.~Tillquist, R.~M.~Frongillo, M.~E.~Lladser,
  Getting the lay of the land in discrete space: A survey of metric dimension and its applications,
  arXiv:2104.07201 [math.CO] (15 Apr 2021).

\bibitem{xu-2021}
  K.~Xu, S.~Klav\v{z}ar,
  Constructing new families of transmission irregular graphs,
  Discrete Appl.\ Math.\ 289 (2021) 383--391.

\end{thebibliography}
\end{document}